\documentclass[twoside]{article}
\usepackage{graphicx}
\usepackage{amsmath}
\usepackage{amssymb}
\newtheorem{theorema}{Theorem.}

\newtheorem{rk}[theorema]{Remark.}
\newcommand\bib[1]{\bibitem[#1]{#1}}

\newcommand\h{h_{\text{\rm top}}}
\newcommand\hm{h_{\text{\rm mult}}}

\begin{document}

\title{A piece-wise affine contracting map\\ with positive entropy}
\author{B. Kruglikov  \& M. Rypdal \\ ~ \\
{\small Institute of Mathematics and Statistics}\\
{\small University of Troms\o, N-9037 Troms\o, Norway}\\
{\small Boris.Kruglikov@matnat.uit.no;
Martin.Rypdal@matnat.uit.no} }
\date{}
\maketitle

\begin{abstract}
We construct the simplest chaotic system with a two-point
attractor.\!\!
\end{abstract}

If $f:X\to X$ is an isometry of the metric space $(X,d)$, then the
topological entropy vanishes: $\h(f)=0$ (for definitions and
notations consult e.g. \cite{4}).

This follows from the fact, that the iterated distance
$d_n^f=\max\limits_{0\le i<n}(f^i)^*(d)$ equals $d$. If $f$ is
distance non-increasing, the same equality holds and again
$\h(f)=0$.

Whenever $f$ can have discontinuities of some tame nature, so that
$f$ is piece-wise continuous, even the isometry result becomes
difficult. In dimension 2 for invertible maps it was proven by
Gutkin and Haydn \cite{3}. In arbitrary dimension Buzzi proved
that piece-wise affine isometries have zero topological entropy
\cite{2}.

In the same paper after the theorem (remark 4) it is claimed that
the result holds for arbitrary piece-wise (non-strictly)
contracting maps. This latter is however wrong and the goal of
this note is to present a counter-example.

\noindent {\bf Example:} Let $X$ be a rhombus $ADBC$ with vertices
$(\pm1,0),(0,\pm1)$, see the figure below. Let $O$ be its center
and $P,Q,R,S$ be on the sides as is shown.

 
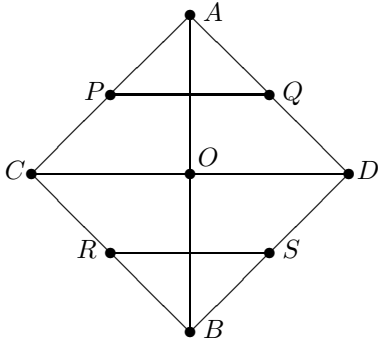
\begin{center}
  \begin{picture}(120,120)
 \put(60,60){\circle*{4}}
 \put(0,60){\circle*{4}}
 \put(120,60){\circle*{4}}
 \put(60,0){\circle*{4}}
 \put(60,120){\circle*{4}}
 \put(30,90){\circle*{4}}
 \put(90,90){\circle*{4}}
 \put(30,30){\circle*{4}}
 \put(90,30){\circle*{4}}
 \put(0,60){\line(1,0){120}}
 \put(60,0){\line(0,1){120}}
 \put(30,90){\line(1,0){60}}
 \put(30,30){\line(1,0){60}}
 \put(0,60){\line(1,1){60}}
 \put(0,60){\line(1,-1){60}}
 \put(120,60){\line(-1,1){60}}
 \put(120,60){\line(-1,-1){60}}
 \put(65,118){$A$}
 \put(65,-2){$B$}
 \put(-10,58){$C$}
 \put(123,58){$D$}
 \put(20,88){$P$}
 \put(95,88){$Q$}
 \put(17,28){$R$}
 \put(95,28){$S$}
 \put(63,63){$O$}
  \end{picture}
 \end{center}

Let $f$ be partially defined on the interior of four big triangles
forming the rhombus. These triangles are bijectively mapped by $f$
as follows:
 $$
ACO\longrightarrow APQ,\ \ ADO\longrightarrow BRS,\ \
BCO\longrightarrow AQP,\ \ BDO\longrightarrow BSR.
 $$
Thus the piece-wise affine map is defined.

If $P,Q$, $R,S$ are middle-points of the intervals $AC,AD$ and
$BC,BD$, then the map is not strictly contracting. But if they are
closer to the vertices $A$ and $B$ respectively than to $C,D$,
then $f$ is strictly contracting. In any case, the attractor of
the system is the 2-point set $\{A,B\}$. Notice that the points
belong to the singularity set, where the map $f$ is not (uniquely)
defined.

Taking $\varepsilon=\frac12$ we observe that the cardinality of
minimal $(n,\varepsilon)$-spanning set satisfies: $2^{n+2}\le
N(f,n,\varepsilon)\le 2^{n+3}$. In fact, if we partition $CD$ into
$2^n$ equal intervals $Z_iZ_{i+1}$, then every $d_n^f$
$\varepsilon$-ball is contained in some triangle $AZ_iZ_{i+1}$ or
$BZ_iZ_{i+1}$ and every such a triangle is covered by two $d_n^f$
$\varepsilon$-balls.

Therefore the topological entropy $\h(f)=\log2$ is positive. In
addition, the Lyapunov spectrum is strictly negative at each point
(for strict contractions), no invariant measure exists
and so the variational principle breaks. \vspace{5pt}

The result of Buzzi \cite{2} generalizes however in the following
fashion:

 \begin{theorema}
Let $f$ be a piece-wise affine map with restriction to each
continuity component being conformal (non-strict) contraction.
Then $\h(f)=0$.
 \end{theorema}

Now we can repeat Buzzi's remark 4 \cite{2}: The proof of his
theorem 3 applies almost literally to the above case of piece-wise
affine conformal contracting maps. Therefore we omit the proof.

 \begin{rk}
It is obvious that if the attractor consists of one point only,
then $\h(f)=0$. If the phase space $X\subset\mathbb{R}^1$ is
one-dimensional and the map is (non-strictly) contracting, then
again $\h(f)=0$. We don't even need to require piece-wise affine
property. This follows from the Buzzi proposition 4 \cite{1},
yielding $\h(f)\le\hm(f)$, where $\hm(f)$ is the multiplicity
entropy, because the latter always vanishes in dimension one.
 \end{rk}
Thus our example with 2 points attractor and 2-dimensional
phase-space $X$ is the simplest possible example with positive
topological entropy.


\end{document}